\documentclass[12pt]{article}%
\usepackage{amsfonts}
\usepackage{sw20bams}
\usepackage{amsmath}
\usepackage{amssymb}
\usepackage{graphicx}%
\providecommand{\U}[1]{\protect\rule{.1in}{.1in}}
\begin{document}

\title{An Exceptional Convolutional Recurrence}
\author{Steven Finch}
\date{September 11, 2024}
\maketitle

\begin{abstract}
A quadratic recurrence of Faltung type, arising via ancestral path lengths of
random binary trees, turns out to be related to the Painlev\'{e} I
differential equation.

\end{abstract}

\footnotetext{Copyright \copyright \ 2024 by Steven R. Finch. All rights
reserved.}Let $T$ be an ordered (strongly) binary tree with $N=2n+1$ vertices.
The \textbf{distance} $\delta(v,w)$ between two vertices $v$ and $w$ of $T$ is
the number of edges in the shortest path connecting them. Let $o$ denote the
root of $T$. The \textbf{height} of a vertex $v$ is $\delta(v,o)$. \ Let
$v\wedge w$ denote the last common ancestor of $v$ and $w$, i.e., the branch
point where the paths from $o$ to $v$ and $w$ diverge. \ Define the
\textbf{internal path length} $p(T)$ and \textbf{ancestral path length} $q(T)$
to be
\[%
\begin{array}
[c]{ccc}%
p(T)=%
%TCIMACRO{\dsum \limits_{v\in T}}%
%BeginExpansion
{\displaystyle\sum\limits_{v\in T}}
%EndExpansion
\delta(v,o), &  & q(T)=%
%TCIMACRO{\dsum \limits_{v,w\in T}}%
%BeginExpansion
{\displaystyle\sum\limits_{v,w\in T}}
%EndExpansion
\delta(v\wedge w,o)
\end{array}
\]
i.e., a sum of $N$, $N^{2}$ heights respectively. \ Janson \cite{Ja-except}
proved that%
\[
\left(  \frac{p(T)}{N^{3/2}},\frac{q(T)}{N^{5/2}}\right)  \rightarrow\left(
\xi,\eta\right)
\]
as $N\rightarrow\infty$, where $\xi$ and $\eta$\ are random variables.
\ Underlying the joint moment
\[
\operatorname*{E}\left(  \xi^{k}\eta^{l}\right)  =\frac{k!l!\sqrt{\pi}%
}{2^{(5k+7l-4)/2}\Gamma((3k+5l-1)/2)}a_{k,l}%
\]
is the following interesting recurrence:
\[
a_{k,l}=2(3k+5l-4)a_{k-1,l}+2(3k+5l-6)(3k+5l-4)a_{k,l-1}+%
%TCIMACRO{\dsum \limits_{0<i+j<k+l}}%
%BeginExpansion
{\displaystyle\sum\limits_{0<i+j<k+l}}
%EndExpansion
a_{i,j}a_{k-i,l-j}%
\]
with $a_{0,0}=-1/2$, $a_{1,0}=1=a_{0,1}$ and $a_{k,l}=0$ when $k<0$ or $l<0$.
All $a_{k,l}$ but $a_{0,0}$ are positive integers when $k\geq0$ and $l\geq0$.
\ We have asymptotics \cite{Ja-except, JC-except}
\[%
\begin{array}
[c]{ccc}%
a_{k,0}\sim\dfrac{1}{2\pi}6^{k}(k-1)!, &  & a_{0,l}\sim C\cdot50^{l}\left(
(l-1)!\right)  ^{2}%
\end{array}
\]
and the identification of the constant $C$ will occupy us immediately.

At the time \cite{F1-except} appeared, high-precision numerical\ analysis by
Kot\v{e}\v{s}ovec \cite{S1-except} gave
\[
C=\frac{\sqrt{15}}{20\pi^{2}}=0.0196207628...=\frac{1}{50}%
(0.9810381421...)=\frac{0.0695537933...}{2\sqrt{\pi}}%
\]
although a rigorous proof still seemed out of reach. \ Setting $k=0=i$ in the
recurrence:%
\[
a_{0,\ell}=2(5l-6)(5l-4)a_{0,l-1}+%
%TCIMACRO{\dsum \limits_{j=1}^{\ell-1}}%
%BeginExpansion
{\displaystyle\sum\limits_{j=1}^{\ell-1}}
%EndExpansion
a_{0,j}a_{0,l-j}%
\]
we see $a_{0,2}=2\cdot4\cdot6\,a_{0,1}+a_{0,1}^{2}=49$. \ Setting $n=\ell-1$:%
\[
a_{0,n+1}=2(5n-1)(5n+1)a_{0,n}+%
%TCIMACRO{\dsum \limits_{j=1}^{n}}%
%BeginExpansion
{\displaystyle\sum\limits_{j=1}^{n}}
%EndExpansion
a_{0,j}a_{0,n+1-j}%
\]
we define $\alpha_{n}=50^{-n}a_{0,n}$. \ It follows that $\alpha_{1}=1/50$,
$\alpha_{2}=49/2500$,%
\[
50^{n+1}\alpha_{n+1}=2\left(  25n^{2}-1\right)  50^{n}\alpha_{n}+%
%TCIMACRO{\dsum \limits_{j=1}^{n}}%
%BeginExpansion
{\displaystyle\sum\limits_{j=1}^{n}}
%EndExpansion
50^{j}\alpha_{j}50^{n+1-j}\alpha_{n+1-j}%
\]
hence%
\begin{align*}
\alpha_{n+1}  &  =\frac{25n^{2}-1}{25}\alpha_{n}+%
%TCIMACRO{\dsum \limits_{j=1}^{n}}%
%BeginExpansion
{\displaystyle\sum\limits_{j=1}^{n}}
%EndExpansion
\alpha_{j}\alpha_{n+1-j}\\
&  =\left(  n^{2}-\frac{1}{25}\right)  \alpha_{n}+2\alpha_{1}\alpha_{n}+%
%TCIMACRO{\dsum \limits_{j=2}^{n-1}}%
%BeginExpansion
{\displaystyle\sum\limits_{j=2}^{n-1}}
%EndExpansion
\alpha_{j}\alpha_{n+1-j}\\
&  =n^{2}\alpha_{n}+%
%TCIMACRO{\dsum \limits_{j=2}^{n-1}}%
%BeginExpansion
{\displaystyle\sum\limits_{j=2}^{n-1}}
%EndExpansion
\alpha_{j}\alpha_{n+1-j}%
\end{align*}
for $n\geq2$. \ Lemma 4 of Joshi \&\ Kitaev \cite{JK-except} confirms that
\[
\lim_{n\rightarrow\infty}\frac{\alpha_{n}}{((n-1)!)^{2}}=c(\alpha_{1})
\]
exists, is finite and is nonzero; not just for $\alpha_{1}=1/50$, $\alpha
_{2}=49/2500$ but more generally for $0<\alpha_{1}<1$ and $\alpha_{2}%
=\alpha_{1}(1-\alpha_{1})$. \ Proposition 18 of \cite{JK-except} further
verifies that $c(1/50)=C$, as expected. \ The paper is devoted to Painlev\'{e}
I. \ A subsequent paper \cite{HJK-except} clarifies that the formula was
originally derived by Takei \cite{Ta-except} and independently by Kitaev.
\ The authors consider Proposition 18 as a straightforward modification of
earlier work \cite{AK-except} on Painlev\'{e} V, thus in a sense it did not
require a separate (lengthy) proof.

The connection to Painlev\'{e} I deserves explanation:\ we provide this in
Section 1, starting with the generating function%
\[
u(t)=%
%TCIMACRO{\dsum \limits_{k=1}^{\infty}}%
%BeginExpansion
{\displaystyle\sum\limits_{k=1}^{\infty}}
%EndExpansion
\alpha_{k}t^{-k}%
\]
and implementing a change of variables. \ A diversion is offered in Section 2:
we discuss Wright's formula \cite{W1-except} for $f(n,n+m)$, the number of
connected graphs with $n$ labeled vertices and $n+m$ edges, and its relation
to the sequence $a_{\cdot,0}$. \ Graphs in this context possess no multiple
edges nor loops.

Another interpretation of the sequences $a_{\cdot,0}$ and $a_{0,\cdot}$
involves Brownian excursion $E(t)$ on the unit interval, which is constrained
to be zero at the endpoints and positive elsewhere. \ The formulas%
\[%
\begin{array}
[c]{ccc}%
\xi=2%
%TCIMACRO{\dint \limits_{0}^{1}}%
%BeginExpansion
{\displaystyle\int\limits_{0}^{1}}
%EndExpansion
E(t)\,dt, &  & \eta=4%
%TCIMACRO{\dint \limits_{0}^{1}}%
%BeginExpansion
{\displaystyle\int\limits_{0}^{1}}
%EndExpansion%
%TCIMACRO{\dint \limits_{0}^{t}}%
%BeginExpansion
{\displaystyle\int\limits_{0}^{t}}
%EndExpansion
\min\limits_{r\leq s\leq t}E(s)\,dr\,dt
\end{array}
\]
guided our treatment in \cite{F1-except}. \ In retrospect, choosing the path
lengths $p$, $q$ would have been more natural. \ 

Sequence $a_{0,\cdot}$ further appears in relation to enumerating rooted maps
(graph drawings) with $n$ edges on an orientable surface of genus $g$.
\ Graphs in this context are allowed to contain multiple edges and loops,
unlike earlier. \ More discussion of this topic is in Section 3. \ For now,
starting with a recurrence \cite{GLM-except}%
\[%
\begin{array}
[c]{lllll}%
b_{0}=1, &  & b_{g+1}=\dfrac{25g^{2}-1}{48}b_{g}-\dfrac{1}{2}%
%TCIMACRO{\dsum \limits_{m=1}^{g}}%
%BeginExpansion
{\displaystyle\sum\limits_{m=1}^{g}}
%EndExpansion
b_{m}b_{g+1-m} &  & \text{for }g\geq0
\end{array}
\]
we see that $b_{1}=-1/48$ and $b_{2}=b_{1}(1-b_{1})/2=-49/(2\cdot48^{2})$.
\ Defining%
\[
\beta_{g}=-\frac{1}{2}\left(  \frac{25}{48}\right)  ^{-g}b_{g}%
\]
it follows that $\beta_{1}=1/50$, $\beta_{2}=49/2500$,%
\begin{align*}
-2\left(  \frac{25}{48}\right)  ^{g+1}\beta_{g+1}  &  =(-2)\left(
\dfrac{25g^{2}-1}{48}\right)  \left(  \frac{25}{48}\right)  ^{g}\beta_{g}-\\
&  \left(  \frac{1}{2}\right)  (-2)^{2}%
%TCIMACRO{\dsum \limits_{m=1}^{g}}%
%BeginExpansion
{\displaystyle\sum\limits_{m=1}^{g}}
%EndExpansion
\left(  \frac{25}{48}\right)  ^{m}\beta_{m}\left(  \frac{25}{48}\right)
^{g+1-m}\beta_{g+1-m}%
\end{align*}
hence%
\begin{align*}
\beta_{g+1}  &  =\left(  \frac{48}{25}\right)  \left(  \dfrac{25g^{2}-1}%
{48}\right)  \beta_{g}+%
%TCIMACRO{\dsum \limits_{m=1}^{g}}%
%BeginExpansion
{\displaystyle\sum\limits_{m=1}^{g}}
%EndExpansion
\beta_{m}\beta_{g+1-m}\\
&  =g^{2}\beta_{g}+%
%TCIMACRO{\dsum \limits_{m=2}^{g-1}}%
%BeginExpansion
{\displaystyle\sum\limits_{m=2}^{g-1}}
%EndExpansion
\beta_{m}\beta_{g+1-m}%
\end{align*}
for $g\geq2$. \ The $\alpha$ and $\beta$ sequences are identical! \ When
writing\ \cite{F2-except}, we failed to notice/report this coincidence. \ For
more occurrences of $a_{0,\cdot}$, with possible variations, see
\cite{Z1-except, Z2-except, FK-except, BOGRW-except, GIKM-except, Zh-except}.

A\ calculation by Kot\v{e}\v{s}ovec \cite{Ko-except} gives%
\[
c(1/2)=0.25338404774059570625310093001580807965353268866793...
\]
but no exact formula for this quantity is known. \ The status of $c(1/50)$ is
exceptional. \ An unrelated problem, also involving the identification of a
constant, is \cite{CFHM-except}
\[
w_{n}\sim6(3.14085756720293695160...)^{-n-1}n
\]
as $n\rightarrow\infty$, where
\[%
\begin{array}
[c]{lllll}%
w_{0}=1, &  & w_{n+1}=\dfrac{1}{(n+1)^{2}}%
%TCIMACRO{\dsum \limits_{k=0}^{n}}%
%BeginExpansion
{\displaystyle\sum\limits_{k=0}^{n}}
%EndExpansion
w_{k}w_{n-k} &  & \text{for }n\geq0.
\end{array}
\]
A solution remains out of reach, as far as is known.

\section{Painlev\'{e}}

Following \cite{HJK-except}, we show first that%
\[
t^{2}\frac{d^{2}u}{dt^{2}}+t\,\frac{du}{dt}-(t+2\alpha_{1})u+t\,u^{2}%
+\alpha_{1}=0.
\]
From%
\[%
\begin{array}
[c]{ccc}%
t\,\dfrac{du}{dt}=-%
%TCIMACRO{\dsum \limits_{k=1}^{\infty}}%
%BeginExpansion
{\displaystyle\sum\limits_{k=1}^{\infty}}
%EndExpansion
k\,\alpha_{k}t^{-k}, &  & t^{2}\dfrac{d^{2}u}{dt^{2}}=%
%TCIMACRO{\dsum \limits_{k=1}^{\infty}}%
%BeginExpansion
{\displaystyle\sum\limits_{k=1}^{\infty}}
%EndExpansion
k(k+1)\alpha_{k}t^{-k},
\end{array}
\]%
\[
t\,u=%
%TCIMACRO{\dsum \limits_{k=1}^{\infty}}%
%BeginExpansion
{\displaystyle\sum\limits_{k=1}^{\infty}}
%EndExpansion
\alpha_{k}t^{-k+1}=%
%TCIMACRO{\dsum \limits_{k=0}^{\infty}}%
%BeginExpansion
{\displaystyle\sum\limits_{k=0}^{\infty}}
%EndExpansion
\alpha_{k+1}t^{-k},
\]%
\[
t\,u^{2}=%
%TCIMACRO{\dsum \limits_{i=1}^{\infty}}%
%BeginExpansion
{\displaystyle\sum\limits_{i=1}^{\infty}}
%EndExpansion%
%TCIMACRO{\dsum \limits_{j=1}^{\infty}}%
%BeginExpansion
{\displaystyle\sum\limits_{j=1}^{\infty}}
%EndExpansion
\alpha_{i}\alpha_{j}t^{-i-j+1}=%
%TCIMACRO{\dsum \limits_{k=1}^{\infty}}%
%BeginExpansion
{\displaystyle\sum\limits_{k=1}^{\infty}}
%EndExpansion%
%TCIMACRO{\dsum \limits_{j=1}^{k}}%
%BeginExpansion
{\displaystyle\sum\limits_{j=1}^{k}}
%EndExpansion
\alpha_{j}\alpha_{k+1-j}t^{-k}%
\]
the coefficient of $t^{-k}$\ on the ODE\ left-hand side becomes%
\[
k^{2}\alpha_{k}-\alpha_{k+1}-2\alpha_{1}\alpha_{k}+%
%TCIMACRO{\dsum \limits_{j=1}^{k}}%
%BeginExpansion
{\displaystyle\sum\limits_{j=1}^{k}}
%EndExpansion
\alpha_{j}\alpha_{k+1-j}=k^{2}\alpha_{k}-\alpha_{k+1}+%
%TCIMACRO{\dsum \limits_{j=2}^{k-1}}%
%BeginExpansion
{\displaystyle\sum\limits_{j=2}^{k-1}}
%EndExpansion
\alpha_{j}\alpha_{k+1-j}=0.
\]
Introducing a change of variables:%
\[%
\begin{array}
[c]{ccc}%
u=\dfrac{1}{2}-\dfrac{1}{2}\left(  \dfrac{x}{6}\right)  ^{-\frac{1}{2}}y, &  &
t=\dfrac{8\sqrt{6}}{25}x^{\frac{5}{2}}%
\end{array}
\]
we compute
\[
\frac{du}{dx}=\frac{1}{24}\left(  \dfrac{x}{6}\right)  ^{-\frac{3}{2}}%
y-\dfrac{1}{2}\left(  \dfrac{x}{6}\right)  ^{-\frac{1}{2}}\frac{dy}{dx},
\]%
\[%
\begin{array}
[c]{ccccc}%
\dfrac{dt}{dx}=\dfrac{4\sqrt{6}}{5}x^{\frac{3}{2}} &  & \text{hence} &  &
\dfrac{dx}{dt}=\dfrac{5\sqrt{6}}{24}x^{-\frac{3}{2}},
\end{array}
\]%
\[
\frac{d^{2}u}{dx^{2}}=-\frac{1}{96}\left(  \dfrac{x}{6}\right)  ^{-\frac{5}%
{2}}y+\frac{1}{12}\left(  \dfrac{x}{6}\right)  ^{-\frac{3}{2}}\frac{dy}%
{dx}-\dfrac{1}{2}\left(  \dfrac{x}{6}\right)  ^{-\frac{1}{2}}\frac{d^{2}%
y}{dx^{2}},
\]%
\[
\dfrac{d^{2}x}{dt^{2}}=-\dfrac{5\sqrt{6}}{16}x^{-\frac{5}{2}}\dfrac{dx}%
{dt}=-\frac{25}{64}x^{-4}.
\]
Substituting%
\[%
\begin{array}
[c]{ccc}%
\dfrac{du}{dt}=\dfrac{du}{dx}\dfrac{dx}{dt}, &  & \dfrac{d^{2}u}{dt^{2}%
}=\dfrac{d^{2}u}{dx^{2}}\left(  \dfrac{dx}{dt}\right)  ^{2}+\dfrac{du}%
{dx}\dfrac{d^{2}x}{dt^{2}}%
\end{array}
\]
into the differential equation for $u(t)$, we obtain%
\[
\dfrac{d^{2}y}{dx^{2}}=6y^{2}+\frac{(50\alpha_{1}-1)y}{4x^{2}}-x.
\]
What is the importance of $\alpha_{1}=1/50$? It is the value for which
Painlev\'{e} I\ emerges.

Using the definition of $u$:%
\[
\dfrac{1}{2}-\dfrac{1}{2}\left(  \dfrac{x}{6}\right)  ^{-\frac{1}{2}}y=%
%TCIMACRO{\dsum \limits_{k=1}^{\infty}}%
%BeginExpansion
{\displaystyle\sum\limits_{k=1}^{\infty}}
%EndExpansion
\alpha_{k}\left(  \dfrac{8\sqrt{6}}{25}x^{\frac{5}{2}}\right)  ^{-k}%
\]
we have%
\[
\dfrac{1}{2}\left(  \dfrac{x}{6}\right)  ^{-\frac{1}{2}}y=-\alpha_{0}-%
%TCIMACRO{\dsum \limits_{k=1}^{\infty}}%
%BeginExpansion
{\displaystyle\sum\limits_{k=1}^{\infty}}
%EndExpansion
\alpha_{k}\left(  \dfrac{8\sqrt{6}}{25}x^{\frac{5}{2}}\right)  ^{-k}%
\]
thus
\[
y=-2\left(  \dfrac{x}{6}\right)  ^{\frac{1}{2}}%
%TCIMACRO{\dsum \limits_{k=0}^{\infty}}%
%BeginExpansion
{\displaystyle\sum\limits_{k=0}^{\infty}}
%EndExpansion
\alpha_{k}\left(  \dfrac{8\sqrt{6}}{25}x^{\frac{5}{2}}\right)  ^{-k}=-\frac
{2}{\sqrt{6}}%
%TCIMACRO{\dsum \limits_{k=0}^{\infty}}%
%BeginExpansion
{\displaystyle\sum\limits_{k=0}^{\infty}}
%EndExpansion
\alpha_{k}\left(  \frac{25}{8\sqrt{6}}\right)  ^{k}x^{-\frac{5k-1}{2}}.
\]
A formal solution $y$ as such is the asymptotic expansion of an actual
Painlev\'{e} I solution (or its analytic continuation, since the series is
everywhere divergent).

\section{Wright}

This diversion involving the sequence $a_{\cdot,0}$ is concerned with the
asymptotics of \cite{W1-except, W2-except, S2-except}%
\[
f(n,n)=\frac{1}{2}\left[  n^{-1}h(n)-n^{n-2}(n-1)\right]  ,
\]%
\[
f(n,n+1)=\frac{1}{24}\left[  n^{n-2}(n-1)\left(  5n^{2}+3n+2\right)
-14h(n)\right]  ,
\]%
\[
f(n,n+2)=\frac{1}{1152}\left[  \left(  45n^{2}+386n+312\right)  h(n)-4n^{n-2}%
(n-1)\left(  55n^{3}+36n^{2}+18n+12\right)  \right]
\]
where%
\[
h(n)=%
%TCIMACRO{\dsum \limits_{k=1}^{n-1}}%
%BeginExpansion
{\displaystyle\sum\limits_{k=1}^{n-1}}
%EndExpansion
\dbinom{n}{k}k^{k}(n-k)^{n-k}.
\]
Clearly $f(3,3)=1=f(4,6)$, as there is exactly one connected labeled graph
with $3$ vertices \& $3$ edges (triangle) and exactly one connected labeled
graph with $4$ vertices \& $6$ edges (tetrahedron). \ Figures 1 and 2 show why
$f(4,4)=15$ and $f(4,5)=6$. \ 

Defining%
\[
\rho_{k-1}=\frac{1}{k!}\operatorname*{E}\left(  \left(  \frac{\xi}{2}\right)
^{k}\right)  =\frac{\sqrt{\pi}}{2^{(7k-4)/2}\Gamma((3k-1)/2)}a_{k,0}%
\]
we have%
\[
\lim_{n\rightarrow\infty}\frac{f(n,n+m)}{n^{n+(3m-1)/2}}=\rho_{m},
\]
i.e., $a_{m+1,0}$ plays a role in approximating the leading coefficient of
graphical counts, for large $n$ and fixed $m$. \ As examples, $a_{1,0}=1$,
$a_{2,0}=5$, $a_{3,0}=60$:%
\[%
\begin{array}
[c]{ccccc}%
\rho_{0}=\dfrac{\sqrt{\pi}}{2^{3/2}\Gamma(1)}a_{1,0}=\dfrac{1}{2}\sqrt
{\dfrac{\pi}{2}}, &  & \rho_{1}=\dfrac{\sqrt{\pi}}{2^{5}\Gamma(5/2)}%
a_{2,0}=\dfrac{5}{24}, &  & \rho_{2}=\dfrac{\sqrt{\pi}}{2^{17/2}\Gamma
(4)}a_{3,0}=\dfrac{5}{128}\sqrt{\dfrac{\pi}{2}}.
\end{array}
\]
Elements of the sequence $\rho_{\cdot}$ are known as \textbf{Wright's
constants} \cite{Jn-except}.%
%TCIMACRO{\FRAME{ftbpFU}{5.3584in}{7.4149in}{0pt}{\Qcb{There are exactly 15
%connected labeled graphs with 4 vertices \& 4 edges.}}{}{except1.eps}%
%width 5.3584in;  height 7.4149in;  depth 0pt;  original-width 8.4207in;
%original-height 11.6949in;  cropleft "0";  croptop "1";  cropright "1";
%cropbottom "0";  filename 'except1.eps';file-properties "NPEU";}}}%
%BeginExpansion
\begin{figure}
[ptb]
\begin{center}
\includegraphics[
height=7.4149in,
width=5.3584in
]%
{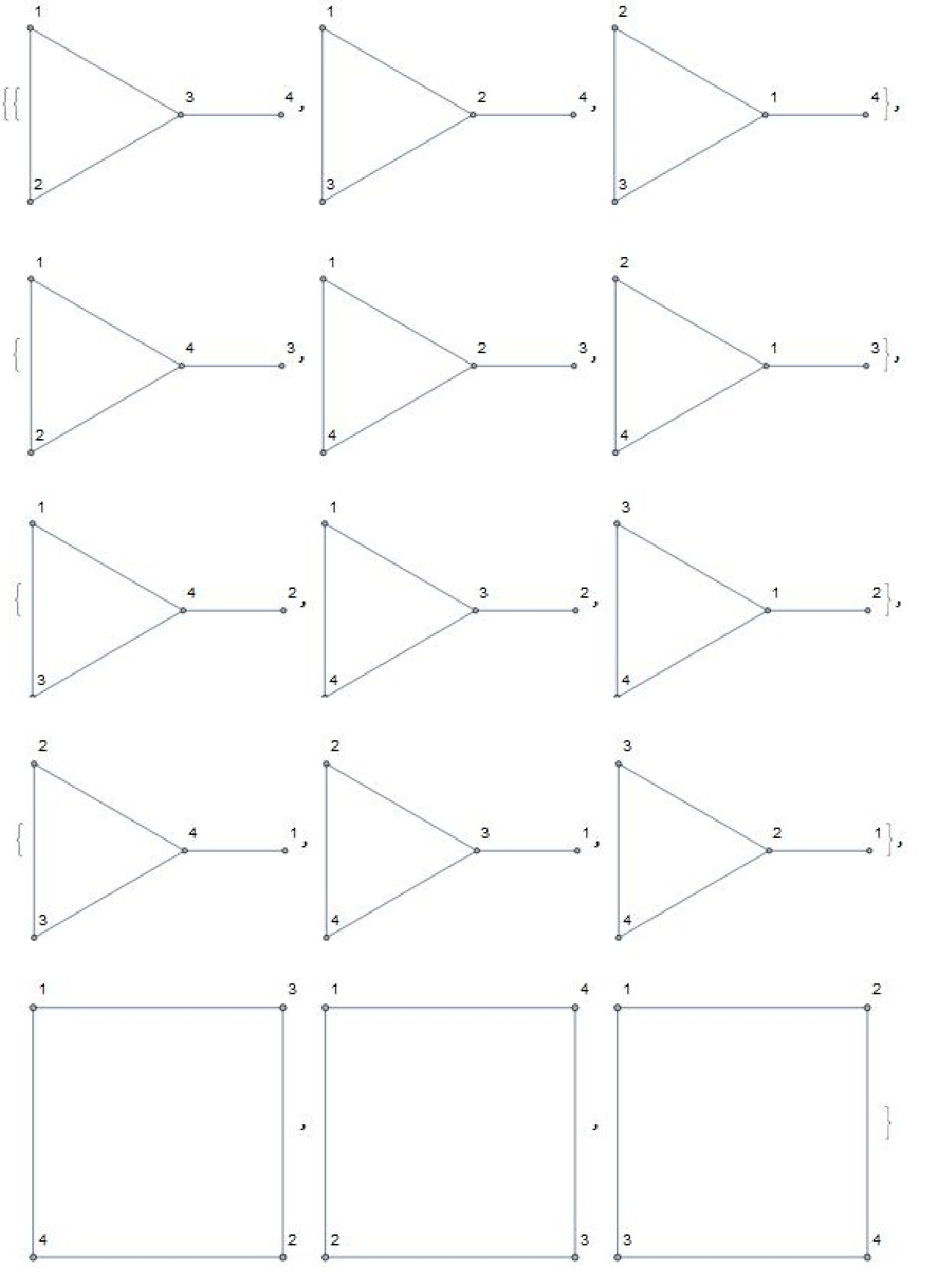}%
\caption{There are exactly 15 connected labeled graphs with 4 vertices \& 4
edges.}%
\end{center}
\end{figure}
%EndExpansion%
%TCIMACRO{\FRAME{ftbpFU}{5.201in}{3.6487in}{0pt}{\Qcb{There are exactly 6
%connected labeled graphs with 4 vertices \& 5 edges.}}{}{except2.eps}%
%width 5.201in;  height 3.6487in;  depth 0pt;  original-width 8.348in;
%original-height 5.8392in;  cropleft "0";  croptop "1";  cropright "1";
%cropbottom "0";  filename 'except2.eps';file-properties "NPEU";}}}%
%BeginExpansion
\begin{figure}
[ptb]
\begin{center}
\includegraphics[
height=3.6487in,
width=5.201in
]%
{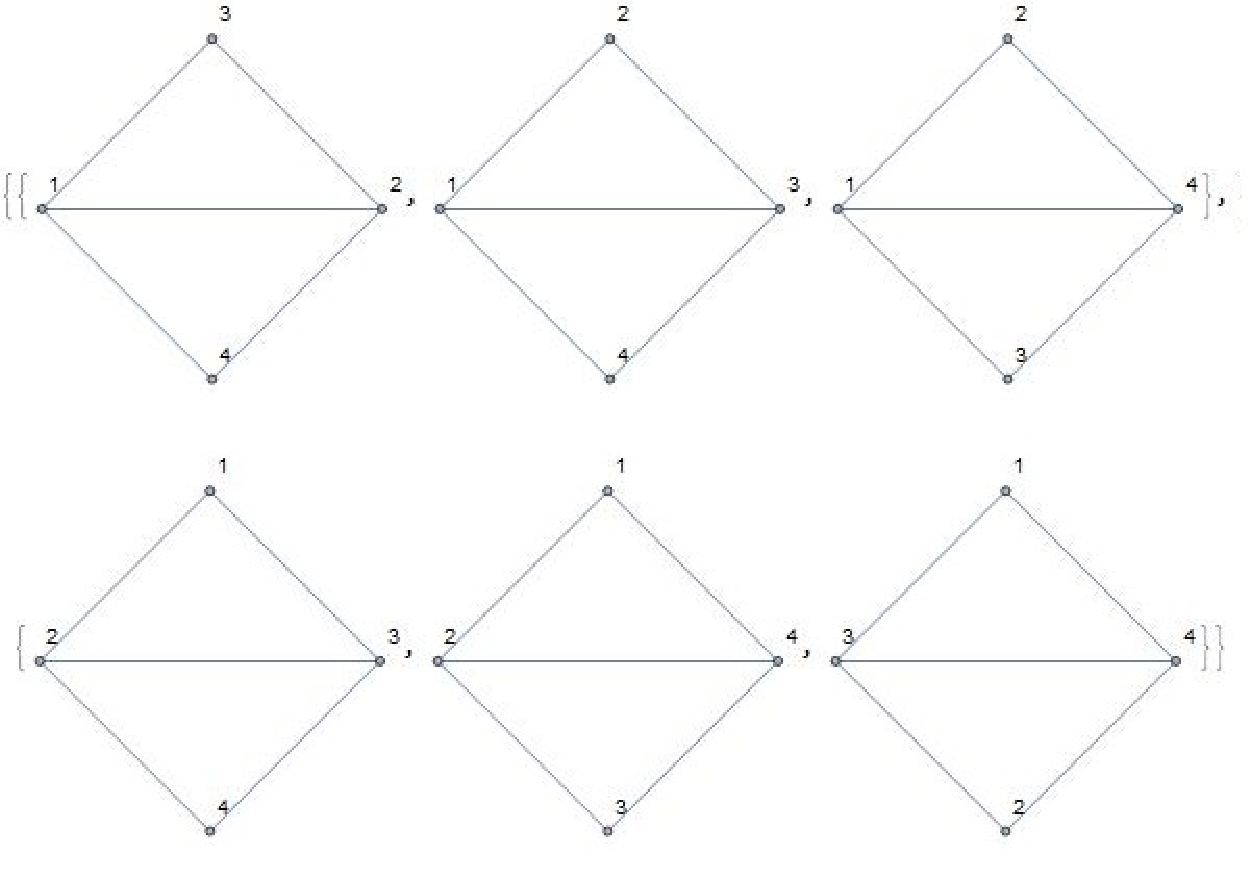}%
\caption{There are exactly 6 connected labeled graphs with 4 vertices \& 5
edges.}%
\end{center}
\end{figure}
%EndExpansion

\section{Tutte}

A \textbf{map} on a compact surface $S$ without boundary is an embedding of a
graph $G$ into $S$ such that all components of $S-G$ are simply connected
\cite{Cf-except}. These components are thus homeomorphic to open disks and are
called \textbf{faces}. The graph $G$ is allowed to have both loops and
multiple parallel edges (unlike those in Section 2). A map is \textbf{rooted}
when an edge, a direction along that edge, and a side of the edge, are
distinguished. The edge is called the \textbf{root edge}, and the face on the
distinguished side is the \textbf{root face}. Two rooted maps are
\textbf{equivalent} if there is a homeomorphism between the underlying
surfaces that preserves all graph incidences and rootedness.

In the case when $S$ is orientable, two rooted maps are equivalent if and only
they are related by an orientation-preserving homeomorphism that (merely)
preserves all graph incidences.\ We shall focus on orientable surfaces only.
The \textbf{genus} $g$ for such surfaces is $0$ for the sphere, $1$ for the
torus, $2$ for the connected sum of two tori, etc. \ 

Let $T_{g}(n)$ denote the number of rooted maps with $n$ edges on an
orientable surface of genus $g$. \ ($T$ stands for "torus" or, equally likely,
Tutte.) \ Figures 3 and 4 show why $T_{0}(1)=2$ and $T_{0}(2)=9$. \ It is
known that $T_{0}(n)$ is the coefficient of $x^{n}$ in the Maclaurin series
expansion \cite{Cf-except, T1-except, T2-except}%
\begin{align*}
\frac{4(2r+1)}{3(r+1)^{2}}  &  =1+2x+9x^{2}+54x^{3}+378x^{4}+2916x^{5}%
+24057x^{6}\\
&  +208494x^{7}+1876446x^{8}+17399772x^{9}+165297834x^{10}+\cdots,
\end{align*}
$T_{1}(n)$ is the coefficient of $x^{n}$ in the expansion
\begin{align*}
\frac{(r-1)^{2}}{12r^{2}(r+2)}  &  =x^{2}+20x^{3}+307x^{4}+4280x^{5}%
+56914x^{6}\\
&  +736568x^{7}+9370183x^{8}+117822512x^{9}+1469283166x^{10}+\cdots,
\end{align*}
$T_{2}(n)$ is the coefficient of $x^{n}$ in the expansion
\begin{align*}
\frac{(r-1)^{4}(r+1)^{2}R}{2304r^{7}(r+2)^{4}}  &  =21x^{4}+966x^{5}%
+27954x^{6}\\
&  +650076x^{7}+13271982x^{8}+248371380x^{9}+4366441128x^{10}+\cdots
\end{align*}
where $r=\sqrt{1-12x}$ and $R=49r^{4}+122r^{3}+225r^{2}+248r+112$. \ Moreover
\cite{BC-except},%
\[
\lim_{n\rightarrow\infty}\frac{T_{g}(n)}{n^{5(g-1)/2}12^{n}}=\tau_{g},
\]
as\thinspace$n\rightarrow\infty$, where
\[
\tau_{g}=\frac{-1}{2^{g-2}\Gamma\left(  (5g-1)/2\right)  }b_{g}%
\]
i.e., $b_{g}$ plays a role in approximating the leading coefficient of rooted
map counts, for large $n$ and fixed $g$. \ As examples, $b_{0}=1$,
$b_{1}=-1/48$, $b_{2}=-49/4608$:%
\[%
\begin{array}
[c]{lllll}%
\tau_{0}=\dfrac{-1}{2^{-2}\Gamma(-1/2)}b_{0}=\dfrac{2}{\sqrt{\pi}}, &  &
\tau_{1}=\dfrac{-1}{2^{-1}\Gamma(2)}b_{1}=\dfrac{1}{24}, &  & \tau_{2}%
=\dfrac{-1}{2^{0}\Gamma(9/2)}b_{2}=\dfrac{7}{4320\sqrt{\pi}}.
\end{array}
\]
Elements of the sequence $\tau_{\cdot}$ are known as \textbf{orientable map
asymptotics constants}. \ A similar sequence can be found for the
non-orientable case, although some aspects of the (difficult) theory remain
conjectural. \
%TCIMACRO{\FRAME{ftbpFU}{5.3246in}{0.768in}{0pt}{\Qcb{There are exactly 2
%rooted maps with 1 edge on a sphere; from \cite{T1-except}.}}{}{except3.eps}%
%width 5.3246in;  height 0.768in;  depth 0pt;  original-width 12.321in;
%original-height 1.7374in;  cropleft "0";  croptop "1";  cropright "1";
%cropbottom "0";  filename 'except3.eps';file-properties "NPEU";}}}%
%BeginExpansion
\begin{figure}
[ptb]
\begin{center}
\includegraphics[
height=0.768in,
width=5.3246in
]%
{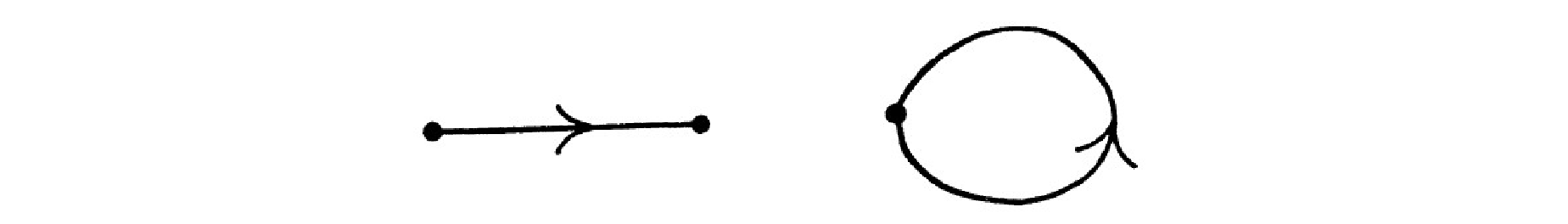}%
\caption{There are exactly 2 rooted maps with 1 edge on a sphere; from
\cite{T1-except}.}%
\end{center}
\end{figure}
%EndExpansion%
%TCIMACRO{\FRAME{ftbpFU}{5.3705in}{2.8496in}{0pt}{\Qcb{There are exactly 9
%rooted maps with 2 edges on a sphere; from \cite{T1-except}.}}{}%
%{except4.eps}{\special{ language "Scientific Word";  type "GRAPHIC";
%width 5.3705in;  height 2.8496in;  depth 0pt;  original-width 12.7543in;
%original-height 6.7369in;  cropleft "0";  croptop "1";  cropright "1";
%cropbottom "0";  filename 'except4.eps';file-properties "NPEU";}}}%
%BeginExpansion
\begin{figure}
[ptb]
\begin{center}
\includegraphics[
height=2.8496in,
width=5.3705in
]%
{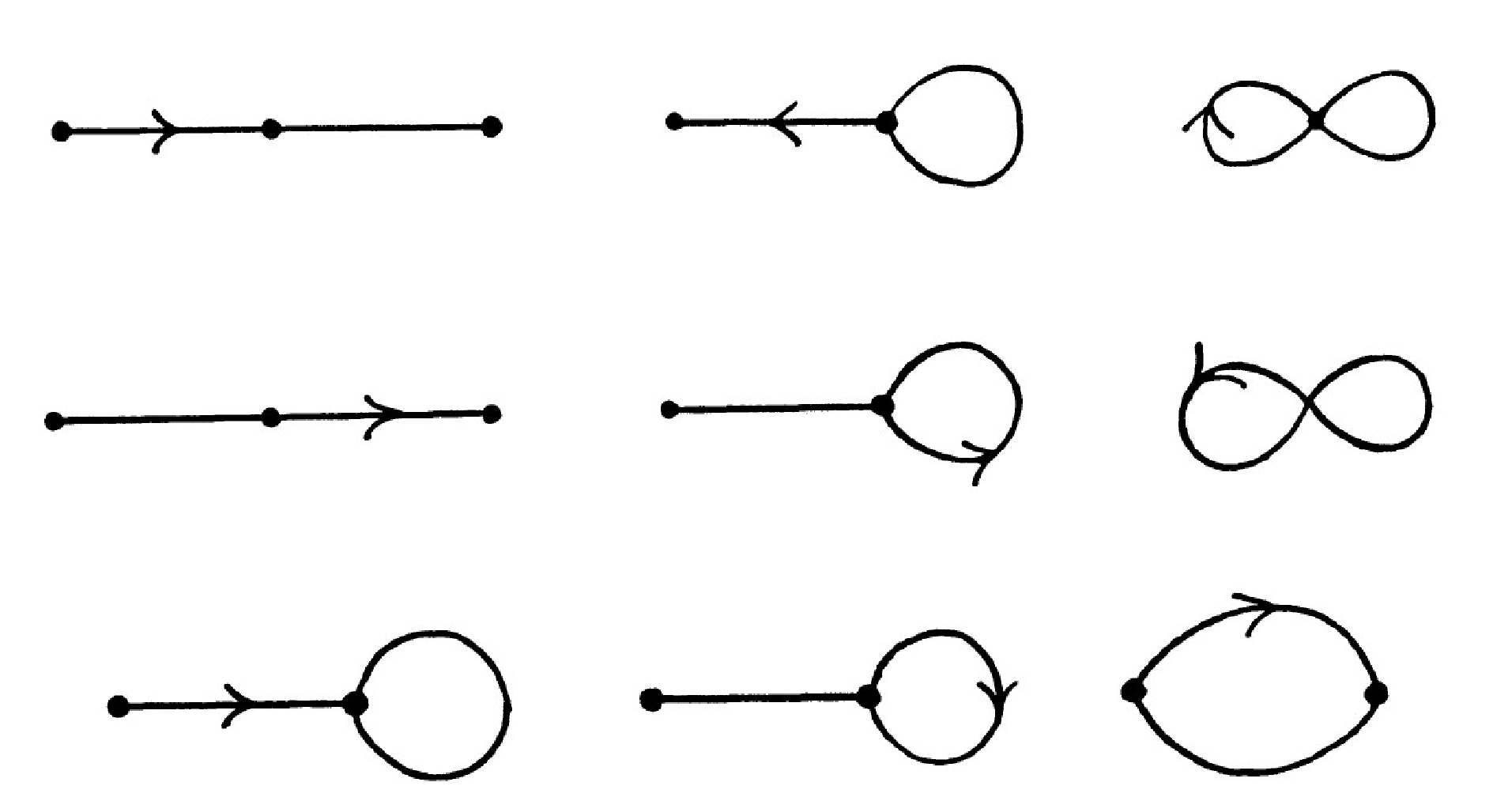}%
\caption{There are exactly 9 rooted maps with 2 edges on a sphere; from
\cite{T1-except}.}%
\end{center}
\end{figure}
%EndExpansion%
%TCIMACRO{\FRAME{ftbpFU}{5.7043in}{1.7123in}{0pt}{\Qcb{(a) An embedding of a
%graph on the torus that is not a map (because the dark face is not simply
%connected). (b) An embedding of the same that is a valid map; from Chapuy
%\cite{Ch-except}.}}{}{except5.eps}{\special{ language "Scientific Word";
%type "GRAPHIC";  maintain-aspect-ratio TRUE;  display "USEDEF";
%valid_file "F";  width 5.7043in;  height 1.7123in;  depth 0pt;
%original-width 7.1382in;  original-height 2.1223in;  cropleft "0";
%filename 'except5.eps';file-properties "NPEU";}}}%
%BeginExpansion
\begin{figure}
[ptb]
\begin{center}
\includegraphics[
height=1.7123in,
width=5.7043in
]%
{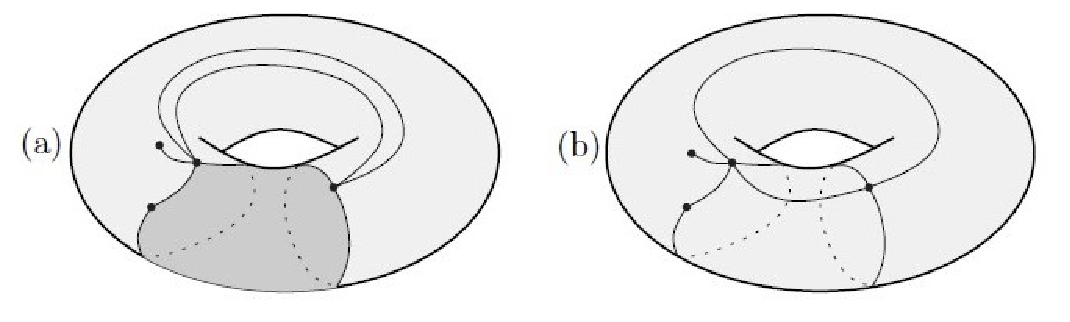}%
\caption{(a) An embedding of a graph on the torus that is not a map (because
the dark face is not simply connected). (b) An embedding of the same that is a
valid map; from Chapuy \cite{Ch-except}.}%
\end{center}
\end{figure}
%EndExpansion%
%TCIMACRO{\FRAME{ftbpFU}{5.3973in}{5.7346in}{0pt}{\Qcb{The reason
%$T_{1}(2)=1\neq2$ is that the two maps are equivalent on the torus (even
%though the indigo loop and the red loop are not homotopic). \ A catalog of all
%$T_{1}(3)=20$ inequivalent maps would be good to see someday.}}{}%
%{except6.eps}{\special{ language "Scientific Word";  type "GRAPHIC";
%width 5.3973in;  height 5.7346in;  depth 0pt;  original-width 14.3758in;
%original-height 15.283in;  cropleft "0";  croptop "1";  cropright "1";
%cropbottom "0";  filename 'except6.eps';file-properties "NPEU";}}}%
%BeginExpansion
\begin{figure}
[ptb]
\begin{center}
\includegraphics[
height=5.7346in,
width=5.3973in
]%
{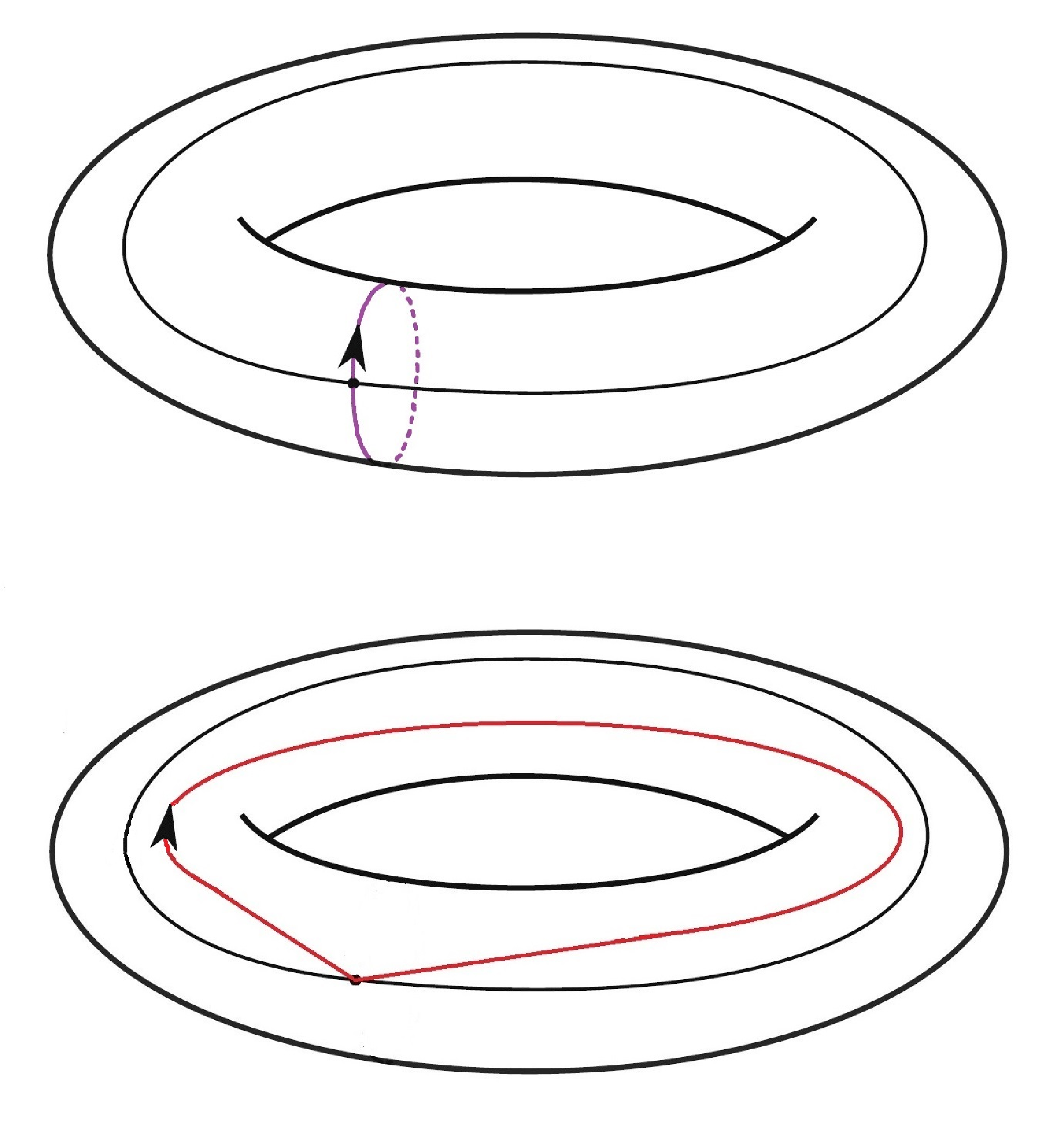}%
\caption{The reason $T_{1}(2)=1\neq2$ is that the two maps are equivalent on
the torus (even though the indigo loop and the red loop are not homotopic).
\ A catalog of all $T_{1}(3)=20$ inequivalent maps would be good to see
someday.}%
\end{center}
\end{figure}
%EndExpansion

\section{Closing Words}

More recent papers \cite{X1-except, X2-except, X3-except} give asymptotics of
the form%
\[
\alpha_{n}\sim\frac{1}{2\pi}\sqrt{\frac{6}{5\pi}}\left(  \frac{1}{4}\right)
^{n}\Gamma\left(  2n-\frac{1}{2}\right)
\]
as $n\rightarrow\infty$ and, for example,%
\[
b_{n}\sim-\frac{1}{\pi}\sqrt{\frac{6}{5\pi}}\left(  \frac{25}{192}\right)
^{n}\Gamma\left(  2n-\frac{1}{2}\right)  ,
\]%
\[%
\begin{array}
[c]{ccc}%
\tau_{n}\sim\dfrac{5}{\pi}\sqrt{\dfrac{6}{5\pi}}\left(  \dfrac{e}%
{1440n}\right)  ^{n/2} &  & \text{(simplifying expressions in
\cite{GLM-except, F2-except}).}%
\end{array}
\]
The phrases \textquotedblleft resurgent methods\textquotedblright\ \&
\textquotedblleft transseries analysis\textquotedblright\ would seem to
suggest a new \&\ active flowering area of research, important for
combinatorialists to assimilate.

In addition to exhibiting the surprising overlap between my two essays
\cite{F1-except} \& \cite{F2-except}, I\ wished to enter into a little more
detail (about Wright's constants and Tutte's drawings). \ The price of writing
with encyclopediac brevity is that, years later, the prose upon re-reading
seems so abstract \& hurried...

\section{Acknowledgements}

I\ am grateful to V\'{a}clav Kot\v{e}\v{s}ovec for his numerical calculation
of $c(1/2)$ and to Alexander Kitaev for a very helpful email discussion.

\end{document}